\documentclass[12pt, draftcls, onecolumn]{IEEEtran}

\ifCLASSINFOpdf
 \else
\fi
\usepackage{cite}
\usepackage{amsmath}
\usepackage{amsthm}
\usepackage{amssymb}
\usepackage{enumitem}
\usepackage{tabularx,booktabs}
\usepackage{multirow}
\usepackage{epsfig}
\usepackage{algorithmic}
\usepackage{algorithm}
\usepackage{xcolor}
\usepackage{array}
\usepackage{xcolor}
\usepackage{multirow,bigdelim}
\usepackage{lipsum}
\usepackage{cite}
\usepackage{amsmath,amssymb,amsfonts}
\usepackage{algorithmic}
\usepackage{graphicx}
\usepackage{textcomp}
\newcolumntype{C}{>{\centering\arraybackslash}X} 
\newtheorem{problem}{Problem}
\newtheorem{thm}{Theorem}
\newtheorem{defn}{Definition}

\newtheorem{lem}{Lemma}

\newtheorem{prop}{Proposition}

\makeatletter
\newcommand\fs@norules{\def\@fs@cfont{\bfseries}\let\@fs@capt\floatc@ruled
  \def\@fs@pre{}%
  \def\@fs@post{}%
  \def\@fs@mid{\kern3pt}%
  \let\@fs@iftopcapt\iftrue}
\makeatother
\restylefloat{algorithm}
\makeatletter
\floatstyle{norules}
\restylefloat{algorithm}

\begin{document}

\title{Characterization of Minimum Time-Fuel Optimal
Control for LTI Systems}

\author{Rajasree~Sarkar,
        Deepak~U.~Patil
        and~Indra~Narayan~Kar
\thanks{Authors are with the Department of Electrical
Engineering, Indian Institute of Technology Delhi, New Delhi, 110016, India.}
\thanks{e-mail: \{Rajasree.Sarkar, deepakpatil, ink\}@ee.iitd.ac.in.}
}

\maketitle

\begin{abstract}
    A problem of computing time-fuel $\int_0^{t_f} (k + |u(t)|)dt$ optimal control for state transfer of a single input linear time invariant (LTI) system to the origin is considered. The input is assumed to be bounded ($|u(t)|\leq 1$). Since, the optimal control is bang-off-bang in nature, it is characterized by sequences of $+1$, $0$ and $-1$ and the corresponding switching time instants. All (candidate) sequences satisfying the Pontryagin's maximum principle (PMP) necessary conditions are characterized. The number of candidate sequences is obtained as a function of the order of system and a method to list all candidate sequences is derived. Corresponding to each candidate sequence, switching time instants are computed by solving a static optimization problem. Since the candidate control input is a piece-wise constant function, the time-fuel cost functional is converted to a linear function in switching time instants. By using a simple substitution of variables, reachability constraints are converted to polynomial equations and inequalities. Such a static optimization problem can be solved separately for each candidate sequence. Finally, the optimal control input is obtained from candidate sequences which gives the least cost. For each sequence, optimization problem can be solved by converting it to a generalized moment problem (GMP) and then solving a hierachical sequence of semidefinite relaxations to approximate the minima and minimizer \cite{Lasserre2008}. Lastly, a numerical example is presented for demonstration of method.
\end{abstract}
\section{Introduction}
Optimal utilization of resources for performing any control task necessitates maximizing the off-duration of control. As a result, problem of computing sparse control for state transfer has gathered a lot of attention over the past few years. Sparse control is of interest in variety of domains, namely, networked, multi-agent control system \cite{bongini2017sparse,LinBop2017sparsenetcon}, transportation systems \cite{Liu2003transport}, etc. One way to obtain sparse control is by computing a control input  that achieves the required state transfer with the least $L_0$-norm (non-zero control duration) \cite{nagahara2016maximum,chatterjee2016characterization}. Such an optimal control problem is difficult to solve because of non-convex and discontinuous nature of the cost function. However, in \cite{nagahara2016maximum}, under certain normality conditions, an equivalence is established  between the solution of $L_0$ and $L_1$ norm optimal control problems. In \cite{challapalli2017continuous}, the authors pose a combination of $L_1$ and $L_2$ cost for getting continuous control inputs along with computational benefits. Further, in \cite{ikeda2019time}, a \emph{time-$L_0$-norm} optimal control problem has also been shown to be equivalent to \emph{time-$L_1$-norm}  (i.e., $\int_0^{t_f} (k + |u(t)|)dt$) optimal control problem. The benefit of such equivalence is that   \emph{time-$L_1$-norm} optimal control problem is tractable compared to the  \emph{time-$L_0$-norm} counterpart. The \emph{time-$L_1$-norm} optimal control problem has received a lot of attention in past \cite{athans1963minimum,athans1964fuel,athans1964fuelrestricted} and is well-known in literature as \emph{time-fuel} optimal control problem \cite{epryan1980, ryan1978synthesis}. 
 
Both time optimal and time-fuel optimal control problem has undergone extensive research in last decade \cite{pontryagin1987mathematical,hajek1971geometric,patil2014computation,athans2013optimal}. However, in comparison to time optimal control problem, computation of open loop as well feedback time-fuel optimal control is a challenging task. Unlike, the time optimal control \cite{pontryagin1987mathematical,patil2014computation}, analytical solutions for time-fuel optimal control cannot be obtained directly from classical control methods like Pontryagin's minimum principle (PMP) even for linear systems except double integrator \cite{athans2013optimal} and a class of second order systems \cite{athans1964fuel,athans1964fuelrestricted}. This is primarily because of two reasons, namely, (i) time optimal control is of bang-bang nature exhibiting transitions only between \emph{two constant levels} i.e. $+1$ and $-1$. On the other hand, time-fuel optimal control shows bang-off-bang nature that changes among \emph{three levels} i.e. $+1$, $0$ and $-1$, and (ii) for real eigenvalues the maximum switching for time optimal control is limited to at most $n-1$ switches \cite{hajek1971geometric} whereas the maximum switches for fuel optimal control is atmost $2n$ \cite{hajek19791} . These drawbacks are seen in various versions of the time-fuel optimal control problem of linear systems that considers free final time \cite{athans1963minimum}, bounded time constraints \cite{athans1964fuelrestricted}, and  \cite{epryan1980}, where authors treat only the fuel optimal control problem for a fixed final time.  The work presented in \cite{kleinman1963fuel} achieves closed loop time-fuel optimal control for a second order system. The work in \cite{kleinman1963fuel} uses all possible sequences of $+1$, $0$ and $-1$, that the optimal control for a second order system follows to constructs state-dependent switching rules. An extension to triple integrator appears in \cite{ryan1978synthesis} in which the following quote appears highlighting difficulty of obtaining a switching rule for a general $n$-th order systems.
\begin{quote}
``the
possibility of obtaining similar results for other third- and higher order
systems seems remote and, in this sense, the approach lacks generality.'' 
\end{quote}
As noted even for second and third order example systems the approach to time-fuel optimal feedback control synthesis is varied. For obtaining a time-fuel optimal feedback control, it is neccessary to find a commonnality that is dependent only on the order of system. Hence, a result that characterizes all possible time optimal control functions that satisfy PMP neccessary conditions can be useful. For a general $n$-th order system \cite{hajek19791} has studied analytical properties of the fuel optimal controls and associated reachability sets with \emph{fixed} final time. 

In this article, assuming \emph{free} final time, we utilize the bang-off-bang property and PMP to characterize sequences of $+1$, $0$ and $-1$, which are identified with candidate time-fuel optimal controls for an $n^{th}$ order LTI system. The sequences that meet the criterion laid by the PMP neccesary conditions are counted and listed. Such a characterization of optimal control candidates is advantageous mainly because, only the knowledge of the order of system is sufficient for listing all possible candidate sequences. It is important, mainly because a list of all possible optimal control candidates can be used as a prior knowledge in formulating any control policy for control problems with time-fuel considerations. To the best of our knowledge, a list and count of all possible time-fuel optimal control candidates for a general $n$-th order LTI system is not available. Further, for each candidate sequence a static non-linear program (NLP) to obtain the time instants at which the optimal control input switch between $+1,0,-1$ is also formulated. The desired control input is obtained by solving several NLPs corresponding to each candidate sequence. The candidate sequence that leads to least time-fuel cost gives the optimal control input. By simple substitution of variables, cost functional can be transformed in to a rational function in decision variables and constraints can be represented by polynomial equations and inequalities. Such a NLP with rational cost function and semi-algebraic constraints has been shown to be equivalent to a generalized moment problem (GMP) in  \cite{Lasserre2008,lasserre2010moments}. Further GMP can be solved by constructing a hierarchy of semidefinite programs \cite{Lasserre2008,lasserre2010moments} using a solver Gloptipoly \cite{henrion2003gloptipoly}. However, this method is not scalable to large problems, and other standard non-linear programming solvers like \emph{fmincon}, SNOPT, IPOPT, etc can also be used. But, these other solvers cannot guarantee global optima. 

It is important to note that direct methods like collocation, and discretization also obtain an NLP, but by \emph{approximating} the original optimal control problem \cite{lin2014control,rao2009survey}. Whereas, indirect methods like shooting method work by \emph{numerically solving} an ordinary differential equation obtained from the PMP. These indirect methods are dependent heavily on the choice of initial co-state guess by the user \cite{rao2009survey} and there are no convergence guarantees available for time-fuel optimal control of a general $n$-th order linear time invariant (LTI) systems. In our approach, we \emph{do not use approximation} at any stage of NLP formulation. Thus, in a way, our approach to obtain time-fuel optimal control is \emph{exact}, in contrast to both direct and indirect methods. The only stage where approximation of solutions takes place is when using solver for NLP.  But, note that in solver Gloptipoly \cite{henrion2003gloptipoly}, the GMP problem obtained from the NLP with rational cost and polynomial cost function is solved by relaxing it to a semidefinite program (SDP) of a finite size called relaxation order. Increasing relaxation order improves the approximation and successive solutions to the increasing sequence of SDP relaxations converges asymptotically to globally optimal solutions \cite{lasserre2010moments}. A preliminary version of this article is published in \cite{sarkar2019computation} where similar results were obtained for second order LTI systems.

\section{Mathematical Notation and Preliminaries}\label{sec:notation}

The \emph{$L_{1}$ norm} of a continuous-time measurable function $u(t)$ over the time interval $[0,T]$, is defined as
$\left \| u(t) \right \|_{1}=\int_{0}^{T}\left | u(t) \right | dt.$
For any set $S$, we define the number of elements in $S$ as its \emph{cardinality} denoted as $\mathcal{N}(S)$.
We are interested in a set of finite-length sequences over $\mathcal{Z}_1:=\{-1,0,1\}$.
For example a set $S=\left\{ (1,0,-1), (1,1,0), (-1,1, -1) \right\}$ is a set of sequences of length 3. 
Consider $A$ and $B$ as two finite length sequence set. Then $A$ is said to be a \emph{sub-sequence set} of $B$ if for every element $a\in A$ there exists an element $b\in B$ such that $a$ is a sub-sequence of $b$. Equivalently, $B$ is termed as the \emph{super-sequence set} of $A$.
For any sequence, its \emph{conjugate sequence} is obtained by reversing the sign of each its elements. For a sequence set $S$, a set of conjugates of all sequences in $S$ is called the \emph{conjugate set} of $S$ and is represented as $\bar{S}$.
For example, consider a sequence $a=(-1,0,1)\in S=\left\{ (1,0,-1), (1,1,0), (-1,1, -1) \right\}$. The conjugate sequence of $a$ is $\bar{a}=(1,0,-1)$. Similarly, the conjugate set of $S$ is $\bar{S}=\left\{(-1,0,1),(-1,-1,0),(1,-1,1)\right\}$.


\section{Problem Formulation}\label{sec: Problem Formulation}
Consider a $n^{th}$ order single input LTI system defined as:
\begin{equation}\label{system_dyn}
\begin{matrix}
\dot{\mathbf{x}}(t)=\mathbf{Ax}(t)+\mathbf{B}u(t),~u(t)\in \mathbb{R}
\end{matrix}
\end{equation}
where $\mathbf{x}(t)\in \mathbb{R}^{n}$ represents the state variable of the system and $\mathbf{A} \in \mathbb{R}^{n\times n}$ and $\mathbf{B} \in \mathbb{R}^{n\times 1}$ are the system and input matrices respectively. Assume that the pair $(\mathbf{A},\mathbf{AB})$ is controllable and the eigenvalues of $\mathbf{A}$ are non-zero, real. Further, for purpose of solving a particular optimization problem in section \ref{sec: Formulation of Optimization Problem}, we will need assumption that eigenvalues are distinct and rational i.e. $(\mathbf{\lambda(A)}\in \mathbb{Q}-\left \{ 0 \right \})$. Hence the eigenvalue $\lambda_{i}$ with numerator $n_{i}$ and denominator $d_{i}$ for $i=1,...,n$ are expressed as $\lambda_{i}=n_{i}/d_{i}=c_{i}/l$ where $c_{i}=n_{i}l/d_{i}$, $l=\text{LCM}(d_{1},...,d_{n}$ and $n_{i},d_{i} \in \mathbb{Z}-\left \{ 0 \right \}$.
Note that such an assumption is not restrictive since rational numbers are dense in the set of real numbers and therefore, any real number can be approximated by a rational number upto arbitrary precision. Later, in Section \ref{sec: Formulation of Optimization Problem}, rational eigenvalue assumption helps in converting the optimal control problem into a tractable optimization problem with rational cost function and constraints described by polynomial inequalities. 
Further, without loss of generality, we assume that $A$ is in diagonal form as $\mathbf{A}=\text{diag}(\mathbf{\lambda(A)})$ and $\mathbf{B}=[b_{1},...,b_{n}]^{T}$ with $b_{1},...,b_{n}\neq 0$. Let the input $u(t)$ be constrained as $\left | u(t) \right |\leq 1$. Thus, the set of the admissible controls is 
\begin{align*}
  U = \left\{u:[0,\infty)\rightarrow [-1,1]\mid \text{u is measurable and }\right.\\
  \left. \left | u(t) \right |\leq 1 \text{ almost everywhere }t \in [0,\infty) \right\}
\end{align*}
Our objective is to choose a control $u(t)\in U$ that steers the system (\ref{system_dyn}) from initial condition $\mathbf{x}(0)=\mathbf{x}_{0}$ to the origin i.e. $\mathbf{x}(t_{f})=0$ with least possible $\left \| u(t) \right \|_{1}$ in finite time $t_{f}$.  To meet this objective, it is necessary that the initial condition $\mathbf{x}_{0}$ lies in the \textit{Reachable Set} defined next.
 The \emph{reachable set} $\mathcal{R}_0$ is the set of all initial conditions $\mathbf{x}_{0}\in \mathbb{R}^{n}$ transferable to the origin by using an admissible control,
$\tiny{
\mathcal{R}_0=\left\{\mathbf{x}_{0} \mid \mathbf{x}_{0}=-\int_{0}^{t} e^{-\mathbf{A}\tau}\mathbf{B}u(\tau)d\tau , u(t)\in U\right\}}.
$
The minimization of  $\left \| u(t) \right \|_{1}$ along with the assurance of finite $t_{f}\geq0$ is formally stated as follows:
\begin{problem}[Time Fuel Optimal Control]\label{prb:TFOCP} Find a control $u(t)\in U$ that steers system \eqref{system_dyn} from $\mathbf{x}_{0}\in \mathcal{R}_0 $
to the origin while minimizing $J=\int_{0}^{t_{f}}\left ( k+\left | u(t) \right | \right )dt, \text{ $t_{f}$: free}$ 
where $k>0$ is a weighing parameter to be appropriately chosen. \end{problem}
A larger $k$ places more weight on time compared to $\left \| u(t) \right \|_{1}$ and thus the control input obtained by solving Problem \ref{prb:TFOCP} will be such that the system trajectory reaches the origin faster with an increased $\left \| u(t) \right \|_{1}$. On the other hand, smaller $k$ is expected to give control input such that state trajectory reaches slowly to the origin, but, also a reduced $\left \| u(t) \right \|_{1}$. The choice of $k$ in the cost function determines the trade-off between the fuel consumption and the speed of system response. Also, $k>0$ ensures that the solution to Problem \ref{prb:TFOCP} drives the state trajectory to the origin in finite time. 
\subsection{Solution To Problem \ref{prb:TFOCP}}\label{sec: Solution To Time Fuel Optimal Problem}
The optimal solution to Problem \ref{prb:TFOCP} necessarily satisfies the  conditions of Pontryagin Minimum Principle (PMP). These conditions will be utilized next for characterizing the candidate functions for the optimal control. We state the Pontryagin Minimum Principle (PMP) as follows \cite{pontryagin1987mathematical}:
\begin{thm}[Pontryagin Minimum Principle (PMP)] \label{thm:pmp}
Let $u^{*}(t)$ be the optimal control function that transfers the initial condition $\mathbf{x}_{0}$ to the origin with minimum cost $J$. Let $\mathbf{x}^{*}(t)$ be the trajectory followed by the system (\ref{system_dyn}) on application of $u^{*}(t)$ with $\mathbf{x}^{*}(0)=\mathbf{x}_{0}$ and $\mathbf{x}^{*}(t_f)=\mathbf{0}$. Then $\mathbf{x}^{*}(t)$ and  $u^{*}(t)$ satisfy the following conditions: 
\begin{enumerate}
\item[(a)] $u^{*}(t)\in U$ is such that for each $t\in[0,t_f]$ it minimizes the Hamiltonian $H$, defined as

\small{$H(\mathbf{x}(t),\mathbf{p}(t),u(t),k):=k+\left|u(t)\right|+\mathbf{p}^{T}(t)(\mathbf{Ax}(t)+\mathbf{B}u(t)),$}
\normalsize
where $\mathbf{p}(t)=[p_{1}(t),...,p_{n}(t)]^{T}$ is the costate, 
\item[(b)] corresponding to $u^{*}(t)$ and $\mathbf{x}^{*}(t)$, there exists an associated optimal costate trajectory $\mathbf{p}^{*}(t)$ which solves the cannonical system:
\begin{align}
\dot{\mathbf{x}}^{*}(t)=\frac{\partial H}{\partial \mathbf{p}}(\mathbf{x}^{*}(t),\mathbf{p}^{*}(t),u^{*}(t),k)\label{state_H}\\ 
\dot{\mathbf{p}}^{*}(t)=-\frac{\partial H}{\partial \mathbf{x}}( \mathbf{x}^{*}(t),\mathbf{p}^{*}(t),u^{*}(t),k)\label{costate}
\end{align}
with boundary conditions: $\mathbf{x}^{*}(0)=\mathbf{x}_{0}$ and $\mathbf{x}^{*}(t_f)=\mathbf{0}$ and 
\item[(c)] terminal condition: $ H( \mathbf{x}^{*}(t),\mathbf{p}^{*}(t),u^{*}(t),k)|_{t=t_f}=0$.
\end{enumerate}
\end{thm}
For $\mathbf{(A,AB)}$ being controllable, the optimal control $u^{*}(t)$ that satisfies all the conditions of Theorem \ref{thm:pmp} is given by:
\begin{equation}\label{red_opt_con}
u^{*}(t)=\left\{\begin{matrix}
-\text{sgn }(\psi(t)) & \text{if } \left | \psi(t) \right |>1,\\
0 & \text{if } \left | \psi(t) \right |< 1,\\
\end{matrix}\right.
\end{equation}
where $\psi(t)=\left \langle \mathbf{B},\mathbf{p}^{*}(t) \right \rangle=b_{1}p_{1}(t)+...+b_{n}p_{n}(t)$. With $\mathbf{A}=$ diag$(c_{i}/l)$ and solution of (\ref{costate}), $\psi(t)$ is further expressed as $\psi(t)=b_{1}\pi_{1}e^{\frac{c_{1}}{l}t}+...+b_{n}\pi_{n}e^{\frac{c_{n}}{l}t}$
where $(\pi_{1},...,\pi_{n})$ is the initial condition of $\mathbf{p}^{*}(t)$.
The initial costate values are unconstrained and as a result, we are unable to determine the optimal control function $u^{*}(t)$. However, we note that $psi(t)$ is a linear combination of several exponential terms. Thus, the following lemma from \cite{pontryagin1987mathematical} can be utilized to characterize the the optimal control function candidates by exploiting the number of roots of functions $\psi(t)$, $\psi(t)-1$ and $\psi(t)+1$.
\begin{lem}\label{lem: real roots}
Let $\eta_{1},\eta_{2},...,\eta_{m}$ be distinct real numbers and let $f_{1}(t),f_{2}(t),...,f_{m}(t)$ be polynomials (with real coefficients) of degree $d_{1},d_{2},...,d_{m}$ respectively. Then the function $f_{1}(t)e^{\eta_{1}t}+f_{2}(t)e^{\eta_{2}t}+...+f_{m}(t)e^{\eta_{m}t}$ has at most $d_{1}+d_{2}+...+d_{m}+m-1$ real roots.
\end{lem}
We define the set of real roots of functions $\psi(t)+1$, $\psi(t)$ and $\psi(t)-1$ as $\psi^{-1}(j) := \{t\geq0~|~\psi(t)=j\}$ for $j=-1,0,+1$ respectively.
Lemma \ref{lem: real roots} combined with \eqref{red_opt_con} helps in concluding that $u^{*}(t)$ is necessarily a piecewise constant function with finitely many switches between $+1$, $0$ and $-1$.
\begin{thm}\label{thm:max dicnt} 
The optimal control function $u^*(t)$ that steers states of system (\ref{system_dyn}) from $x_{0} \in \mathcal{R}_0$ to the origin with minimum $J$ satisfies the following conditions (c.f. \cite{hajek19791}):
\begin{enumerate}
    \item [(i)] $u^*(t)$ is a  piecewise constant function on an interval $t \in [0,t_f]$ switching between $+1$, $0$ and $-1$. Moreover, switching always takes place between $+1$ to $0$, $0$ to $+1$, $-1$ to $0$ and $0$ to $-1$. Direct switching between $+1$ to $-1$ and vice-a-versa is not possible. 
    \item [(ii)] $u^{*}(t_f)\neq 0$ and is always equal to $+1$ or $-1$,
    \item [(iii)] The function $\psi(t)$ that defines $u^*(t)$ is such that
    \begin{enumerate}\label{thm:max dicnt3} 
     \item [(a)] $\mathcal{N}(\psi^{-1}(0))\leq n-1$
    \item [(b)] $\mathcal{N}(\psi^{-1}(+1))\leq n$,
    \item [(c)] $\mathcal{N}(\psi^{-1}(-1))\leq n$,
\end{enumerate}
    \item [(iv)] $u^*(t)$ has at most $2n$ discontinuities.
\end{enumerate}
\end{thm}
\begin{proof}
The proof of each statement in the theorem is provided in separate ordered arguments:
\begin{enumerate}
    \item [(i)] From \eqref{red_opt_con}, we see that $u^{*}(t)$ is a piecewise constant function switching between $+1$, $0$ and $-1$ values determined by $\psi(t)$. Also, let $[\tau_0,\tau_1]$ be an interval such that $\psi(\tau_0)=-1$ and  $\psi(\tau_1)=+1$. By continuity of $\psi(t)$, there exists a time $\tau_2\in[\tau_0,\tau_1]$ such that $\psi(\tau_2)=0$. Thus, there exists a sub interval in $[\tau_0,\tau_1]$ in which the $|\psi(t)|\leq 1$ and as a result corresponding input $u^{*}(t)=0$ on that sub-interval.  Thus, the corresponding $u^{*}(t)$ transits between the $-1$ to $+1$ values through zero. Similarly it can be shown that $u^{*}(t)$ transits between the $+1$ to $-1$ values also through zero.
    \item [(ii)] If $u^{*}(t_f)$ is zero, then $ H( \mathbf{x}^{*}(t),\mathbf{p}^{*}(t),u^{*}(t),k)|_{t=t_f}=k$, thus, violating condition (c) in Theorem \ref{thm:pmp}. Therefore, $u^{*}(t_f)\neq 0$. From (i), it follows that $u^{*}(t_f)$ is equal to $\pm1$.
\item [(iii)] From Lemma \ref{lem: real roots}, we note that $\psi(t)=0$ has $n-1$ real roots. Thus, $\mathcal{N}(\psi^{-1}(0))\leq n-1$. Next, let us define
\small{\begin{align*}
\psi_{+1}(t):=\psi(t)-1=b_{1}\pi_{1}e^{\frac{c_{1}}{l}t}+...+b_{n}\pi_{n}e^{\frac{c_{n}}{l}t}-e^{\frac{c_{n+1}}{l}t}\\
\psi_{-1}(t):=\psi(t)+1=b_{1}\pi_{1}e^{\frac{c_{1}}{l}t}+...+b_{n}\pi_{n}e^{\frac{c_{n}}{l}t}+e^{\frac{c_{n+1}}{l}t}
\end{align*}}
\normalsize where $c_{n+1}=0$. Recall that $\lambda(A)$ was assumed to be non-zero distinct rational. Therefore, $c_{1},...,c_{n},c_{n+1}$ are distinct real numbers. Again from Lemma \ref{lem: real roots}, we conclude that $\psi_{+1}(t)$ and $\psi_{-1}(t)$ has $n$ real roots respectively. Thus, $\mathcal{N}(\psi^{-1}(+1))\leq n$ and $\mathcal{N}(\psi^{-1}(-1))\leq n$.   

\item [(iv)] Note from \eqref{red_opt_con} that for the function $u^{*}(t)$  switch happens only when $\psi(t)=+1$ or $\psi(t)=-1$. Hence, by using (iii), the number of discontinuities in $u^{*}(t)$ is $\mathcal{N}(\psi^{-1}(+1))+\mathcal{N}(\psi^{-1}(-1))\leq n+n =2n$. \qed
\end{enumerate}
\end{proof}
Theorem \ref{thm:max dicnt}, gives necessary conditions that optimal control candidates must satisfy.
From Theorem \ref{thm:max dicnt}, the resulting form of $u^{*}(t)$ can be expressed as follows:
\begin{equation}\label{opt_con_time_inst}
\small{    u(t)=\left\{\begin{matrix}
\vdots  & \vdots \\ 
0  & t\in(t_{m},t_{m+1})\\ 
-1  & t\in(t_{m+1},t_{m+2})\\ 
0 & t\in(t_{m+2},t_{m+3})\\ 
+1  & t\in(t_{m+3},t_{m+4})\\ 
\vdots  & \vdots\\
\pm 1 & t\in(t_{f-1},t_{f})
\end{matrix}\right.}
\end{equation} 
where, switching time instances satisfy $
t_{1}<...< t_{m} < t_{m+1}< ... < t_{f}.$ Let $\mathcal{U}\subset U$ be the set of all possible inputs of the form given by \eqref{opt_con_time_inst} i.e., $\mathcal{U}=\{u:[0,t_f]\to \mathcal{Z}_1 \mid u \mbox{ is piecewise constant}\}$. Further, let the set of all inputs of the form \eqref{opt_con_time_inst} satisfying all the conditions of Theorem \ref{thm:max dicnt} be denoted as $\mathcal{U}^* \subset \mathcal{U}$. In other words, $\mathcal{U}^*$ is the set of all control inputs that satisfy the PMP necessary conditions. 
\subsection{Correspondence between input set and sequence set}
Note that if we ignore the values of switching instants in \eqref{opt_con_time_inst} and consider constant values arranged in their temporal order, then each element of $\mathcal{U}$ can be represented by a finite length sequence over $\mathcal{Z}_1=\{-1,0,1\}$.  For example input of the form \eqref{opt_con_time_inst} can be compactly represented by the following sequence $(..., 0, +1, 0, -1, ..., \pm 1)$
This allows us to define an equivalence relation among elements of $\mathcal{U}$ defined as follows: 
\begin{defn}
Two inputs  $u_1(t),u_2(t)\in \mathcal{U}$ are equivalent if their corresponding sequences are the same.
\end{defn}
 This equivalence relation divides $\mathcal{U}$ into disjoint equivalence classes. Moreover, a sequence over $\mathcal{Z}_1$ represents all equivalent inputs belonging to that respective equivalence class. For example  all inputs, \[\small{u(t)=\left\{\begin{matrix}-1, & t\in[0,t_1)\\ 0, & t\in[t_1,t_2)\\ 1 ,& t\in [t_2,t_3)\end{matrix}\right.}\] with switching instants satisfying $0<t_1<t_2<t_3<\infty$ are equivalent and are represented by a sequence $(-1,0,1)$. Consequently,  sequence $(-1,0,1)$ forms an equivalence class of all inputs with different values of the time instants in above mentioned form.  Thus, a bijective map can be set up between the set of equivalence classes of $\mathcal{U}$ and the set of sequences over $\mathcal{Z}_1$.

Since, all conditions of Theorem
\ref{thm:max dicnt} put restrictions only on the number of switching events and temporal order in which $+1$, $-1$ and $0$ appear in optimal control candidates, it is easier to deal with sequences rather than piecewise continuous functions.  
Hence, to characterize all possible optimal control candidates that satisfy conditions of Theorem
\ref{thm:max dicnt}, we list all the equivalence classes to which optimal control candidates belong to. We call the sequence for which the corresponding equivalent inputs satisfy the conditions of Theorem
\ref{thm:max dicnt} as a candidate sequence. 



\begin{table}[ht!]
\caption{Number of times the three levels of $\psi(t)$ are crossed for certain sub-sequence}
\label{tab:CRs}
\begin{center}
\begin{tabular}{cccc}
\hline
{Sub-sequence} &  {$\mathcal{N}(\psi^{-1}(+1))$} & $\min \mathcal{N}(\psi^{-1}(0))$ & {$\mathcal{N}(\psi^{-1}(-1))$}\\
          \hline
$(0, 1)$       & \multicolumn{1}{c}{0} & \multicolumn{1}{c}{0} &\multicolumn{1}{c}{1} \\
$(1, 0)$          & \multicolumn{1}{c}{0} & \multicolumn{1}{c}{0} &\multicolumn{1}{c}{1} \\
$(0, -1)$         & \multicolumn{1}{c}{1} & \multicolumn{1}{c}{0} &\multicolumn{1}{c}{0} \\
$(-1, 0)$          & \multicolumn{1}{c}{1} & \multicolumn{1}{c}{0} &\multicolumn{1}{c}{0} \\
$(1, 0 , -1)$          & \multicolumn{1}{c}{1} & \multicolumn{1}{c}{1} &\multicolumn{1}{c}{1} \\
$(-1,0, +1)$          & \multicolumn{1}{c}{1} & \multicolumn{1}{c}{1} &\multicolumn{1}{c}{1} \\\hline
\end{tabular}
\end{center}
\end{table}

\subsection{Candidate Sequence Properties}\label{sec: Candidate Sequences}
Candidate sequences for optimal control $u^*(t)$ are sequences obtained by arranging $+1$'s, $0$'s and $-1$'s in various combinations that satisfy the conditions of Theorem \ref{thm:max dicnt}. This section describes the structure of such candidate sequences.

\subsubsection{Candidate Sequence Structure}\label{sec: Structural Framework}
Consider the following sub-sequences: (a) $(+1,0)$, (b) $(-1,0)$, (c) $(0,+1, 0)$, (d) $(0, -1, 0)$, (e) $(0, +1)$, (f) $(0, -1)$, (g) $(+1)$ and (h) $(-1)$. A finite concatenation of these sub-sequences yields another sequence. To be consistent with conditions (i), (ii), (iii-b) and (iii-c) of Theorem \ref{thm:max dicnt}, the concatenation of these elements should satisfy the following requirements.
\begin{enumerate}
    \item Terminal sub-sequence should end with a non-zero value,
    \item any two consecutive sub-sequences should be such that, the first sub-sequence ends with zero value and the next sub-sequence starts with zero.
\end{enumerate}
Further, to satisfy condition (iii) of Theorem \ref{thm:max dicnt}, we would subsequently consider the number of times $\psi(t)$ crosses $+1$, $0$ and $-1$ corresponding to each sub-sequence (shown in Table \ref{tab:CRs}). Now, in the resulting sequence obtained from such a concatenation, we combine two consecutive zeroes into a single zero. For example, a concatenation of the sub-sequences that satisfy the above conditions, is done as follows:
\begin{align*}
\small{    (0, +1, \underbrace{0),(0}, +1, \underbrace{0),(0}, -1, \underbrace{0),(0}, -1)}
\end{align*}
The sequence extracted will be $(0, +1, 0, +1, 0, -1, 0,-1).$
Let $\mathcal{S}$ be the set of all sequences obtained by finite concatenation of sub-sequences (a)-(h). Note that any sequence in $\mathcal{S}$ already satisfies condition (i) of Theorem \ref{thm:max dicnt}.

To characterize a candidate sequence we must further identify sequences from $\mathcal{S}$ that satisfy remaining conditions of Theorem \ref{thm:max dicnt}. For that, we decompose a candidate sequence into the following three segments 
\begin{enumerate}
    \item \textit{Beginning segment:}  A candidate sequence begins with any one of the sub-sequences $(+1,0)$, $(-1,0)$, $(0,+1, 0)$ and $(0, -1, 0)$. We denote $\mathcal{N}(\psi^{-1}(+1))$ in this segment as $\beta^{+}$ and $\mathcal{N}(\psi^{-1}(-1))$ as $\beta^{-}$. The various values of $(\beta^{+},\beta^{-})$ for these four beginning sub-sequences are $(1,0)$, $(0,1)$, $(2,0)$ and $(0,2)$ respectively. 
\item \textit{Middle segment:} This segment is a finite concatenation of sub-sequences that both start and end with zero values, namely, $(0, +1, 0)$ and $(0 , -1 ,0)$. Let this segment be a concatenation of $\gamma^{+}$ and $\gamma^{-}$ numbers of $(0, +1, 0)$ and $(0, -1 , 0)$ sub-sequences respectively. Then $\mathcal{N}(\psi^{-1}(-1))$  and $\mathcal{N}(\psi^{-1}(+1))$ in this segment are $2\gamma^{+}$ and $2\gamma^{-}$ respectively.
\item \textit{End segment:} This segment is formed out of sub-sequences $(0, +1)$ and $(0 , -1)$ to satisfy condition (ii) of Theorem \ref{thm:max dicnt}. Let $\mathcal{N}(\psi^{-1}(+1))=\epsilon^{+}$ and $\mathcal{N}(\psi^{-1}(-1))=\epsilon^{-}$. The values of $(\epsilon^{+},\epsilon^{-})$ for these possible transitions are $(1,0)$ and $(0,1)$ respectively.
\end{enumerate}
Let us denote the set of sequences in $\mathcal{S}$ for which $\mathcal{N}(\psi^{-1}(+1))=p$ and $\mathcal{N}(\psi^{-1}(-1))=q$ by $S_{p,q}\subset \mathcal{S}$. Also, let $S^{+}_{p,q}$  be the set of sequences in $S_{p,q}$ that begin with sub-sequences $(+1,0)$ or $(0,+1,0)$ and similarly, $S_{p,q}^{-}$ be the set of sequences that begin with $(-1,0)$ or $(0,-1,0)$ in $S_{p,q}$. Then, $S_{p,q} = S^{+}_{p,q} \cup S^{-}_{p,q}$ and $S^{+}_{p,q}\cap S^{-}_{p,q} = \emptyset$.
Moreover, for any sequence in $S^{+}_{p,q}$, we have $\beta^+\in\{1,2\}$ and $\beta^-=0$. Similarly for  any sequence in $S^{-}_{p,q}$, we have $\beta^-\in\{1,2\}$ and $\beta^+=0$.
Therefore, for any sequence in $S_{p,q}$ with $p+q>1$, the quantities $(\beta^{+},\epsilon^{+},\gamma^{+},\beta^{-},\epsilon^{-},\gamma^{-})$ must satisfy the following set of equations and inequalities.
\begin{subequations} \label{eqn}
  \begin{align}
  & \beta^{+}+2\gamma^{+}+\epsilon^{+} = q,~\gamma^{+}\in \mathbb{Z}_{\geq 0} \label{eqn1}\\
    & \beta^{-}+2\gamma^{-}+\epsilon^{-} = p,~\gamma^{-}\in \mathbb{Z}_{\geq 0} \label{eqn2}\\
    & \beta^{+}\beta^{-}=0, ~\beta^{+},\beta^{-}\in \left\{0,1,2\right\}\label{eqn3}\\
    &  \beta ^{+}+\beta ^{-}\neq 0 \label{eqn5}\\
   & \epsilon^{+}\epsilon^{-} = 0,~\epsilon^{+},\epsilon^{-}\in \left\{0,1\right\}\label{eqn4}\\
    & \epsilon^{+}+\epsilon^{-} = 1 \label{eqn6}
\end{align}
\end{subequations} 
The solution set for equation \eqref{eqn} is listed in Table \ref{tab:crossings} for all cases of $p$ and $q$.
\begin{table}[h!]
\caption{Solution set of system of equations \eqref{eqn}}
\label{tab:crossings}
\begin{center}
\begin{tabular}{ccccccccc}
\hline

Set & $p$ & $q$ & \multicolumn{1}{c}{$\beta^{+}$} & \multicolumn{1}{c}{$\epsilon^{+}$} &\multicolumn{1}{c}{$\gamma^{+}$} & \multicolumn{1}{c}{$\beta^{-}$} &\multicolumn{1}{c}{$\epsilon^{-}$} &\multicolumn{1}{c}{$\gamma^{-}$} \\\hline 

\multirow{4}{*}{$S^{+}_{p,q}$} & \multicolumn{1}{c}{even,$\geq 0$} & \multicolumn{1}{c}{odd,$\geq 3$} &\multicolumn{1}{c}{2} &\multicolumn{1}{c}{1} &\multicolumn{1}{c}{$\frac{q-3}{2}$} &\multicolumn{1}{c}{0} &\multicolumn{1}{c}{0} &\multicolumn{1}{c}{$\frac{p}{2}$} \\
                            
    & \multicolumn{1}{c}{odd,$\geq 1$} & \multicolumn{1}{c}{even,$\geq 2$}  &\multicolumn{1}{c}{2} &\multicolumn{1}{c}{0} &\multicolumn{1}{c}{$\frac{q-2}{2}$} &\multicolumn{1}{c}{0} &\multicolumn{1}{c}{1} &\multicolumn{1}{c}{$\frac{p-1}{2}$}\\ 
    
    & \multicolumn{1}{c}{odd,$\geq 1$} & \multicolumn{1}{c}{odd,$\geq 1$}  &\multicolumn{1}{c}{1} &\multicolumn{1}{c}{0} &\multicolumn{1}{c}{$\frac{q-1}{2}$} &\multicolumn{1}{c}{0} &\multicolumn{1}{c}{1} &\multicolumn{1}{c}{$\frac{p-1}{2}$} \\
    
    & \multicolumn{1}{c}{even,$\geq 0$} & \multicolumn{1}{c}{even,$\geq 2$} &\multicolumn{1}{c}{1} &\multicolumn{1}{c}{1} &\multicolumn{1}{c}{$\frac{q-2}{2}$} &\multicolumn{1}{c}{0} &\multicolumn{1}{c}{0} &\multicolumn{1}{c}{$\frac{p}{2}$} \\ \hline 
    
\multirow{4}{*}{$S^{-}_{p,q}$} & \multicolumn{1}{c}{even,$\geq 2$} & \multicolumn{1}{c}{odd,$\geq 1$} &\multicolumn{1}{c}{0} &\multicolumn{1}{c}{1} &\multicolumn{1}{c}{$\frac{q-1}{2}$} &\multicolumn{1}{c}{2} &\multicolumn{1}{c}{0} &\multicolumn{1}{c}{$\frac{p-2}{2}$} \\
                            
    & \multicolumn{1}{c}{odd,$\geq 3$} & \multicolumn{1}{c}{even,$\geq 0$}  &\multicolumn{1}{c}{0} &\multicolumn{1}{c}{0} &\multicolumn{1}{c}{$\frac{q}{2}$} &\multicolumn{1}{c}{2} &\multicolumn{1}{c}{1} &\multicolumn{1}{c}{$\frac{p-3}{2}$}\\ 
    
    & \multicolumn{1}{c}{odd,$\geq 1$} & \multicolumn{1}{c}{odd,$\geq 1$} &\multicolumn{1}{c}{0} &\multicolumn{1}{c}{1} &\multicolumn{1}{c}{$\frac{q-1}{2}$} &\multicolumn{1}{c}{1} &\multicolumn{1}{c}{0} &\multicolumn{1}{c}{$\frac{p-1}{2}$} \\ 
   
    & \multicolumn{1}{c}{even,$\geq 2$} & \multicolumn{1}{c}{even,$\geq 0$}  &\multicolumn{1}{c}{0} &\multicolumn{1}{c}{0} &\multicolumn{1}{c}{$\frac{q}{2}$} &\multicolumn{1}{c}{1} &\multicolumn{1}{c}{1} &\multicolumn{1}{c}{$\frac{p-2}{2}$} \\ \hline 
\end{tabular}
\end{center}
\end{table}
For the remaining sequences that arise when $p+q\leq1$, we directly use Table \ref{tab:CRs} to obtain the following sequence sets: (i) $S^+_{0,0}=\left\{(+1)\right\}$, (ii) $S^-_{0,0}=\left\{(-1)\right\}$, (iii) $S^+_{0,1}=\left\{(0,+1)\right\}$, (iv) $S^-_{0,1}=\varnothing$, (v) $S^+_{1,0}=\varnothing$, (vi) $S^-_{1,0}=\left\{(0,-1)\right\}$.

Next, we note the symmetry in the set $S_{p,q}$ with $p+q>1$. Flipping the sign of each non-zero element that appear in a sequence in set $S_{p,q}^+$ gives us a sequence in $S_{q,p}^-$. In other words, conjugate set (recall from section \ref{sec:notation}) $\overline{S_{p,q}^+}\subset S_{q,p}^-$. Similarly, each sequence in set $S_{q,p}^-$ is also a sequence in $\overline{S_{p,q}^+}$, giving us $ S_{q,p}^-\subset\overline{S_{p,q}^+}$.  

\begin{prop}\label{lem:subsets_rel}
$S^{-}_{q,p}=\overline{S^{+}_{p,q}},~p+q > 1$
\end{prop}


Note from Proposition \ref{lem:subsets_rel}, that we only need to list sequences in $S^{+}_{p,q}$. Thus, solving the following system of equations to realize the sequence structure in $S^{+}_{p,q}$ for $p+q>1$ suffices to construct the set $S_{p,q}$. 
\begin{subequations} \label{eqn:modified_eqn}
\begin{align}
   & \beta^{+}+2\gamma^{+}+\epsilon^{+} = q,~\gamma^{+}\in \mathbb{Z}_{\geq 0},~\beta^{+}\in \left\{1,2\right\} \label{eqn:modified_eqn1}\\
    & 2\gamma^{-}+\epsilon^{-} = p,~\gamma^{-}\in \mathbb{Z}_{\geq 0} \label{eqn:modified_eqn2}\\
    & \epsilon^{+}\epsilon^{-} = 0,~\epsilon^{+},\epsilon^{-}\in \left\{0,1\right\}\label{eqn:modified_eqn3}\\
     & \epsilon^{+}+\epsilon^{-} = 1 
\end{align}
\end{subequations}
The feasible solution for equations \eqref{eqn:modified_eqn} are shown in Table \ref{tab:crossings}. Recall that each solution of \eqref{eqn:modified_eqn} corresponds  to a particular sequence in set $S_{p,q}^+$.

\subsection{Computation of Number of Candidate Sequences}
In this subsection, we aim to compute the number of candidate sequences. At first, we propose the following lemma that gives a count on the number of sequences in $S_{p,q}$.
\begin{lem}\label{lem:number}
Let $p\leq n$ and $q\leq n$ be fixed. Then, the number of sequences in $S_{p,q}$ for $p+q>1$, is $\mathcal{N}(S_{p,q})=\mathcal{N}(S^{+}_{p,q}) + \mathcal{N}(S^{-}_{p,q})$ with $\mathcal{N}(S^{+}_{p,q})=^{(\gamma^{+}_{p,q}+\gamma^{-}_{p,q})}C_{\gamma^{+}_{p,q}}$ and $\mathcal{N}(S^{-}_{p,q}) = ^{(\gamma^{+}_{q,p}+\gamma^{-}_{q,p})}C_{\gamma^{+}_{q,p}}$.
\end{lem}
\begin{proof}
In a solution to \eqref{eqn:modified_eqn} for a fixed value of $p\leq n$ and $q\leq n$, the values of $\beta^+,\beta^-$ and $\epsilon^+, \epsilon^-$ are also fixed as per Table \ref{tab:crossings}. This means that the starting and end segments are same for all sequences in $S_{p,q}^+$ for that $p$ and $q$. As a consequence, the only distinguishing factor to separate one sequence from another in $S^{+}_{p,q}$ is the arrangement of sub-sequences in their middle segments. Hence, the number of sequences in $S^{+}_{p,q}$ is determined by computing the number of possible concatenations of $\gamma^{+}_{p,q}$ number of $(0, +1, 0)$ sub-sequences  and $\gamma^{-}_{p,q}$ number of  $(0, -1, 0)$ sub-sequences which is $\mathcal{N}(S^{+}_{p,q})=    ^{(\gamma^{+}_{p,q}+\gamma^{-}_{p,q})}C_{\gamma^{+}_{p,q}}$. Similarly, for the set $S_{p,q}^-=\overline{S}^{+}_{q,p}$ (from Proposition \ref{lem:subsets_rel}), the number of possible concatenations of $\gamma^{+}_{q,p}$ $(0,-1,0)$ sub-sequences, and $\gamma^{-}_{q,p}$  $(0,+1, 0)$ sub-sequences is $\mathcal{N}(S^{-}_{p,q})=^{(\gamma^{+}_{q,p}+\gamma^{-}_{q,p})}C_{\gamma^{+}_{q,p}}$. Since $S^{+}_{p,q}\cap S^{-}_{p,q} = \emptyset$, we get $\mathcal{N}(S_{p,q})=\mathcal{N}(S^{+}_{p,q}) + \mathcal{N}(S^{-}_{p,q})$. \qed
\end{proof}

Note that sequences in $S_{p,q}^+$ for any $p,q\leq n$ satisfy condition (i), (ii), (iii-b) and (iii-c) in Theorem \ref{thm:max dicnt}. However, it is to be ascertained whether all sequences also satisfy condition (iii-a) in Theorem \ref{thm:max dicnt}. Note that, since condition (iii-a) restricts number of roots of $\psi(t)$ and does not affect the switching transitions for optimal input, it is impossible to verify condition (iii-a) directly for all $S^+_{p,q}$ with $p+q>1$, except the one sequence $\tilde{u}^+ = (+1, 0, -1, 0, (-1)^2, \hdots,(-1)^{n-1},0, (-1)^{n})$ in $S^+_{n,n}$, which is a concatenation of $n$ sequences $(+1,0,-1)$ with alternating signs. 
Clearly $\tilde{u}^+$ will lead to $\mathcal{N}(\psi^{-1}(0))=n$ and thus violating condition (iii-a) of Theorem \ref{thm:max dicnt}. Similarly, for $S^-_{n,n}$, the sequence $\tilde{u}^-$ can be obtained which violates the condition (iii-a) of Theorem \ref{thm:max dicnt}. For rest of the sequences, the minimum of $\mathcal{N}(\psi^{-1}(0))<n$ can be ensured and hence condition (iii-a) of Theorem \ref{thm:max dicnt} is satisfied.

Let $S^+=\bigcup_{p=0}^{n}\bigcup_{q=0}^{n} S_{p,q}^+$ for $p+q>1$. The set $S^{+}$ can be divided into four disjoint sets, namely, $S^{+}_{1},S^{+}_{2},S^{+}_{3},S^{+}_{4}$ defined as:
\begin{align*}
    S^{+}_{1} &= \left\{s\in S^{+}\mid \beta^{+}=2,\epsilon^{+}=1, \epsilon^{-}=0\right\}\\
    S^{+}_{2} &= \left\{s\in S^{+}\mid \beta^{+}=2,\epsilon^{+}=0, \epsilon^{-}=1\right\}\\
    S^{+}_{3} &= \left\{s\in S^{+}\mid \beta^{+}=1,\epsilon^{+}=0, \epsilon^{-}=1\right\}\\
    S^{+}_{4} &= \left\{s\in S^{+}\mid \beta^{+}=1,\epsilon^{+}=1, \epsilon^{-}=0\right\}
\end{align*}
Thus, we now get following theorem to count all the candidate sequences. 
\begin{thm}
The number of candidate sequences in $S^+$ is
\begin{align*}
   \left\{ \begin{matrix}
^nC_{\frac{n-2}{2}}+2~^nC_{\frac{n}{2}}+~^{n+1}C_{\frac{n}{2}}-3 & \text{for }n \text{ even}\\ 
3~^nC_{\frac{n-1}{2}}+~ ^{n+1}C_{\frac{n+1}{2}}-3 & \text{for }n \text{ odd}
\end{matrix}\right.
\end{align*}
\end{thm}
\begin{proof}
Using Table \ref{tab:crossings} and Lemma \ref{lem:number}, the number of sequences in $S_{1}^{+}$ for odd $n$ is sum of $^{\frac{p+q-3}{2}}C_{\frac{q-3}{2}}$ over $p=0,2,...,n-1$ nested with $q=3,5,...,n$. We simplify this sum by substituting $i=p/2$ and $j=(q+1)/2$ as follows:
\begin{align}\label{seq_count_exp}
\sum_{i=0}^{\frac{n-1}{2}}\sum_{j=2}^{\frac{n+1}{2}}~ ^{i+j-2}C_{i}
\end{align}
Since, $\sum_{k=0}^{p}~^{r+k}C_{k}= ~^{r+p+1}C_{p}$, we get
\begin{align*}
\sum_{j=2}^{\frac{n+1}{2}}~^{j-1+\frac{n-1}{2}}  C_{\frac{n-1}{2}}=\sum_{j=1}^{\frac{n+1}{2}}~^{j-1+\frac{n-1}{2}}  C_{\frac{n-1}{2}}-1
\end{align*}
Therefore, $\sum_{j=0}^{\frac{n-1}{2}}~^{j+\frac{n-1}{2}}C_{j} -1 ~=~ ^{n}C_{\frac{n-1}{2}}-1$.
Following similar mathematical operations, the total number of sequences in $S_i^+,i=2,3,4$ are $^{n}C_{\frac{n-1}{2}}-1$, $^{n+1}C_{\frac{n+1}{2}}-1$ and $^{n}C_{\frac{n-1}{2}}-1$ respectively. Identical manipulations can be made for $n$ even and number of sequences in $S_i^+,i=1,2,3,4$ are $^{n}C_{\frac{n-2}{2}}-1$, $^{n}C_{\frac{n}{2}}-1$, $^{n}C_{\frac{n}{2}}-1$ and $^{n+1}C_{\frac{n-1}{2}}-1$ respectively.
By eliminating $\tilde{u}^{+}$ and considering the sequences in $S_{p,q}^{+}$ with $p+q\leq 1$, the theorem statement follows.\qed
\end{proof}
\section{Characterization of Time-fuel Optimal Control}\label{sec: Generalized Structure of Optimal Control}
Recall the definition of super-sequence set from Section \ref{sec:notation}. The super-sequence set for $S_i^+,i=1,..,4$ sets are summarized in Table \ref{table: basis_seq_even} and \ref{table: basis_seq_odd} for $n$ even and $n$ odd respectively.
\begin{table}[h!]
\caption{Super-sequential Set of Equivalence Classes of $S^{+}$ for $n=$ even}
\label{table: basis_seq_even}
\begin{center}
\begin{tabular}{c c c m {0.6cm} m{0.6cm} c}
\hline 
Equiv. & Super & Starting & \multicolumn{2}{c}{Intermittent Seg.} & End\\
Class   &  Set  & Segment  & $\gamma^{+}$ & $\gamma^{-}$                           & Segment\\
\hline
$S_{1}^{+}$ & $S_{n,n-1}^{+}$ & $(0, +1 , 0)$  & $\frac{n-4}{2}$ & $\frac{n}{2}$ & $(0, +1)$\\
$S_{2}^{+}$ & $S_{n-1,n}^{+}$ & $(0, +1 , 0)$  & $\frac{n-2}{2}$ & $\frac{n-2}{2}$ & $(0, -1)$\\
$S_{3}^{+}$ & $S_{n-1,n-1}^{+}$ & $(+1 , 0)$  & $\frac{n-2}{2}$ & $\frac{n-2}{2}$ & $(0, -1)$\\
$S_{4}^{+}$ & $S_{n,n}^{+}$ & $(+1 , 0)$  & $\frac{n-2}{2}$ & $\frac{n}{2}$ & $(0, +1)$\\
\hline
\end{tabular}
\end{center}
\end{table}
\begin{table}[h!]
\caption{Super-sequential Set of Equivalence Classes of $S^{+}$ for $n=$ odd}
\label{table: basis_seq_odd}
\begin{center}
\begin{tabular}{c c c m {0.6cm} m{0.6cm} c}
\hline 
Equiv. & Super & Starting & \multicolumn{2}{c}{Intermittent Seg.} & End\\
Class   &  Set  & Segment  & $\gamma^{+}$ & $\gamma^{-}$                           & Segment\\
\hline
$S_{1}^{+}$ & $S_{n-1,n}^{+}$ & $(0, +1 , 0)$  & $\frac{n-3}{2}$ & $\frac{n-1}{2}$ & $(0, +1)$\\
$S_{2}^{+}$ & $S_{n,n-1}^{+}$ & $(0, +1 , 0)$  & $\frac{n-3}{2}$ & $\frac{n-1}{2}$ & $(0, -1)$\\
$S_{3}^{+}$ & $S_{n,n}^{+}$ & $(+1 , 0)$  & $\frac{n-1}{2}$ & $\frac{n-1}{2}$ & $(0, -1)$\\
$S_{4}^{+}$ & $S_{n-1,n-1}^{+}$ & $(+1 , 0)$  & $\frac{n-3}{2}$ & $\frac{n-1}{2}$ & $(0, +1)$\\
\hline
\end{tabular}
\end{center}
\end{table}


Thus, the set of all candidate sequences in $S^+$ is $S^+_{n,n}\cup S^+_{n-1,n} \cup S^+_{n,n-1} \cup S^+_{n-1,n-1}-\left\{\tilde{u}^+\right\}$. Note that for $n>2$, the set $S^+_{n-1,n}\cup S^+_{n,n}-\left\{\tilde{u}^+\right\}$ forms the super-sequence set of $S^+_{n,n-1},S^+_{n-1,n-1}$ and therefore also includes all candidate sequences. But, for $n=2$, $S^+_{2,2}-\left\{\tilde{u}^+\right\}=\varnothing$, so,  $S_{0,2}$ is considered to be the super-sequence of $S^+_{4}$ for $n=2$. Therefore, we represent the set of all candidate sequences in $S^+$ as:
\begin{align}\label{opt_set}
   \left\{\begin{matrix}
    S^+_{n-1,n} \cup S^+_{0,n} & \text{for}~n=2 \\
S^+_{n-1,n} \cup S^+_{n,n}-\left\{\tilde{u}^+\right\} & \text{for}~n>2 
\end{matrix}\right.
\end{align}

For $n=2$, we get only two sequences $(0,+1,0,-1)\in S^+_{1,2}$  and $(+1,0,+1)\in S^+_{0,2}$.
 For $n>2$, all elements of the set $S_{n,n}^+-\left\{\tilde{u}^+\right\}$ are the subsequences of the following sequence,
\begin{align}\label{opt_seq1}
     (+1,0,s_1,0,...,s_{n-1},0,(-1)^n)
\end{align}
with $\sum_{m=1}^{n-1}s_m=-1$ (for $n$ even) and $\sum_{m=1}^{n-1}s_m=0$ (for $n$ odd) with $s_{m}^2-1=0$ for $m=1,...,n-1$. Similarly, all elements of the set $S_{n-1,n}^+$  are the subsequences of the following sequence,
\begin{align}\label{opt_seq2}
    (0,+1,0,s_1,0,...,s_{n-2},0,-(-1)^n)
\end{align}
where $s_{m}^2-1=0$  for $m=1,...,n-2$, $\sum_{m=1}^{n-2}s_m=0$ (for $n$ even) and $\sum_{m=1}^{n-2}s_m=-1$ (for $n$ odd). 
We combine \eqref{opt_seq1} and \eqref{opt_seq2} to construct the following  function $u_c^+(t)$ from which all other time-fuel optimal control candidates (with $+1$ appearing as the beginning non-zero input) can be obtained by putting equality constraints on the consecutive time instants.
\begin{align}\label{gen_opt_con}
    & u_c^{+}(t)=\left\{\begin{matrix}
0 & t\in\left [ 0,t_{1} \right ]\\
+1 & t\in\left [ t_{1},t_{2} \right ]\\
0 & t\in\left [ t_{2},t_{3} \right ]\\
s_{1} & t\in\left [ t_{3},t_{4} \right ]\\
\vdots & \vdots \\
0 & t\in\left [ t_{2n-2},t_{2n-1} \right ]\\
s_{n-1} & t\in\left [ t_{2n-1},t_{2n} \right ]\\
0 & t\in\left [ t_{2n},t_{2n+1} \right ]\\
(-1)^{n} & t\in\left [ t_{2n+1},t_{2n+2} \right ]\\
\end{matrix}\right. 
\end{align}
\begin{align}\label{mipconstraints}
    \begin{matrix}
\text{where}, & s_{i}^{2}=1~\forall i=1,\dots,n-1\\
& \sum_{i=1}^{n-1}s_{i}=
\left\{\begin{matrix}
0 & \text{for }n\text{ odd}\\ 
-1 & \text{for }n\text{ even},
\end{matrix}\right.\\
& 0 \leq t_1\leq t_2 \leq...\leq t_{2n+1}\leq t_{2n+2}
\end{matrix}
\end{align}
and at least two pairs of switching time instances $t_j$ for $j=1,...,2n+2$ are equal. Let us denote the set of all inputs of the form \eqref{gen_opt_con} satisfying constraints \eqref{mipconstraints} by $\mathcal{U}^+_c$. Similarly, we  obtain $u_c^-(t)$ and the corresponding set $\mathcal{U}^-_c$. Finally, note that $\mathcal{U}^*\subset \mathcal{U}_c^+\cup\mathcal{U}_c^-$.



\section{Optimal Control Problem to Optimization Problem}\label{sec: Formulation of Optimization Problem}

In this section, we use $\mathcal{U}_c^{\pm}$ obtained in Section \ref{sec: Generalized Structure of Optimal Control} and convert Problem \ref{prb:TFOCP} to an equivalent optimization problem.
\subsection{Constraints}
The set of initial sets that can be steered to the origin using inputs from $\mathcal{U}_c^+$ is $\mathcal{X}^+_c=\{\mathbf{x}_0=-\int_{0}^te^{-At}Bu(t)~dt,~u(t)\in \mathcal{U}_c^+\}.$ Note that since the input is in piecewise constant form, the term $-\int_{0}^te^{-At}Bu(t)~dt$ is a sum of exponentials of parameters $t_{1},...,t_{2n+2}$. Recall, $A$ is in diagonal form with $\lambda(A)$ along the diagonal. Therefore, $\mathcal{X}^+_c$ can be alternately represented in terms of $e^{\lambda_{1}t},...,e^{\lambda_{n}t}$. Using $\lambda_{i} = c_{i}/l$, $i=1,...,n$ and performing a substitution as follows:
\begin{equation}\label{eqn: substitution}
    a_{j}=e^{\frac{t_{j}}{l}} \text{ for } j=1,...,2n+2,
\end{equation}
the parametric representation of $\mathcal{X}^+_c$ is achieved in terms of $a_{j}$'s. This representation is polynomial if all the eigenvalues i.e. $\lambda_{i}$'s have same sign. However, $\lambda_{i}$'s with both positive and negative signs result in a rational parametric representation of $\mathcal{X}^+_c$. Therefore, the expression for each component of $\mathbf{x}_0$, in general, is written as $x_{i,0}=N_{i}(a_{1},...,a_{2n+2})/D_{i}(a_{1},...,a_{2n+2})=N_{i,}/D_{i}$
where $N_{i}(a_{1},...,a_{2n+2})$ is the numerator and  $D_{i}(a_{1},...,a_{2n+2})$ is the denominator polynomial of the rational function $x_{i,0}$ for $i=1,...,n$. Rearranging, a polynomial equality constraint is obtained as $x_{i,0}D_i-N_i = 0$.
Note that substitution \eqref{eqn: substitution} translate the inequality constraint $0\leq t_1 \leq ...\leq t_{2n+2}$ to $1\leq a_1\leq ...\leq a_{2n+2}$.
\subsection{Cost function}
Using the piecewise constant nature of $u(t)\in \mathcal{U}_c^+$, the cost function $J$ in Problem \ref{prb:TFOCP} is expressed as a weighted combination of final time and time duration for which $u(t)$ is non-zero. The cost function with $u(t)$ as \eqref{opt_seq1}, denoted by $J_1$, is $J_1 = kt_{2n+1}+t_1-t_2+t_3....-t_{2n}+t_{2n+1}-1$. Subsequently, the cost function with $u(t)$ as \eqref{opt_seq2}, denoted as $J_2$, is $J_2 = kt_{2n}-t_1+t_2-t_3-....-t_{2n-1}+t_{2n}$
With substitution \eqref{eqn: substitution}, and by monotonically increasing nature of logarithms, the cost function is expressed as
\begin{equation*}\label{without_log_cost}
    J_1= \frac{a_{1}a_{3}\dots a_{2n+1}^{k+1}}{a_{2}a_{4}\dots a_{2n}},~\text{and}~J_2= \frac{a_{2}a_{4}\dots a_{2n}^{k+1}}{a_{1}a_{3}\dots a_{2n-1}} 
\end{equation*}

\subsection{Time-fuel optimization problem}\label{time_fuel_opt_prb}
With the constraint and cost function defined above, we are required to solve two sets of optimization problem for $n>2$ as:
\begin{align*}
\tag{OP1}
   \text{Minimize}~~ & J_1 =\frac{a_{1}a_{3}\dots a_{2n+1}^{k+1}}{a_{2}a_{4}\dots a_{2n}}\nonumber\\
   \text{Subject to}~~ & x_{i,0}=- \frac{b_{i}l}{c_{i}}\left [-1+a_{1}^{-c_{i}}-s_1 a_{2}^{-c_{i}}+s_{1}a_{3}^{-c_{i}}- \right.\nonumber \\
   & ~~~~~~~~~\dots-s_{n-1}a_{2n-2}^{-c_{i}} + s_{n-1}a_{2n-1}^{-c_{i}} -  \nonumber\\
   & ~~~~~~~~~(-1)^{n}a_{2n}^{-c_{i}}+(-1)^{n}a_{2n+1}^{-c_{i}}\left.\right ],~i=1,...,n,\nonumber\\
   & a_{j}-a_{j+1}\leq 0, ~\forall j=1,\dots,2n, \nonumber\\
   & \sum_{m=1}^{n-1}s_{m}=
\left\{\begin{matrix}
-1 & \text{for }n\text{ even}\\ 
0 & \text{for }n\text{ odd},
\end{matrix}\right.\nonumber\\
& a_{1}\geq 1,~s_{m}\in \left\{+1,-1 \right\}
\end{align*}
\begin{align*}
\tag{OP2}
   \text{Minimize}~~ & J_2 =\frac{a_{2}a_{4}\dots a_{2n}^{k+1}}{a_{1}a_{3}\dots a_{2n-1}}\nonumber\\
   \text{Subject to}~~ & x_{i,0}=- \frac{b_{i}l}{c_{i}}\left [-a_{1}^{-c_{i}}+ a_{2}^{-c_{i}}-s_{1}a_{3}^{-c_{i}}+s_{1}a_{4}^{-c_{i}} \right.\nonumber \\
   & ~~~~~~~~~\dots-s_{n-2}a_{2n-3}^{-c_{i}} + s_{n-2}a_{2n-2}^{-c_{i}} + \nonumber\\
   & ~~~~~~~~~(-1)^{n}a_{2n-1}^{-c_{i}}-(-1)^{n}a_{2n}^{-c_{i}}\left.\right ],~i=1,...,n \nonumber\\
   & a_{j}-a_{j+1}\leq 0, ~\forall j=1,\dots,2n-1, \nonumber\\
   & \sum_{m=1}^{n-2}s_{m}=
\left\{\begin{matrix}
0 & \text{for }n\text{ even}\\ 
-1 & \text{for }n\text{ odd},
\end{matrix}\right.\nonumber\\
& a_{1}\geq 1,~s_{m}\in \left\{+1,-1 \right\}
\end{align*}

\subsection{Discussion}
Note that problems (OP1) and (OP2) are mixed-integer nonlinear programming problems (MINLP), which are in general computationally difficult to solve even with the available solvers. Therefore, we treat these MINLP's as a collection of multiple optimization problems by putting the values of the integer variable $s_{m},~m=1,2,...$, in the constraints. Based on the values of $s_m$ and eliminating one corresponding to $\tilde{u}^+$, the number of optimization problems that are required to be solved in the form (OP1) is $^{n-1}C_{\frac{n}{2}}-1$ for $n$ even (and $^{n-1}C_{\frac{n-1}{2}}-1$ for $n$ odd). Similarly, the number of optimization problems that are required to be solved in the form (OP2) is $^{n-2}C_{\frac{n-2}{2}}$ for $n$ even (and $^{n-2}C_{\frac{n-1}{2}}$ for $n$ odd). Note that these two sets of optimization problems are obtained for $u(t)\in \mathcal{U}_c^{+}$. Two more sets of such problems can be formulated for $u(t)\in \mathcal{U}_{c}^-$ in similar manner. Therefore, we are required to solve, in total, the following number of non-linear programs for $n\geq 2$ with rational cost function and semi-algebraic constraints,
\begin{align*}
    2\left (^{n-1}C_{\frac{n}{2}}-1+^{n-2}C_{\frac{n-2}{2}}\right) ~~~& \text{for}~n\text{ even}\\
    2\left (^{n-1}C_{\frac{n-1}{2}}-1+^{n-2}C_{\frac{n-1}{2}}\right) ~~~& \text{for}~n\text{ odd}.
\end{align*}
The time-fuel optimal control is obtained by solving these optimization problems. Each problem can be solved by converting it into a generalized moment problem and approximating it by a hierarchy of  semidefinite programs (See \cite{lasserre2010moments} for more details). After solving all the optimization problems of the form (OP1) and (OP2), there is a possibility of multiple solutions yielding the same minimum cost. All solutions which yields minimum cost are selected and depending on requirements in terms of number of switchings, time of state-transfer and the $L_1$ norm of input, a suitable optimal solution can be chosen. Also, since we are required to solve each optimization problem separately, the computation can done in a distributed manner. We also note finally that the existence of solution for at least one problem is guaranteed if and only if $\mathbf{x}_{0}\in \mathcal{R}_0$. 
\subsubsection{Solver for the optimization problem}
Each optimization problem being defined with a rational cost function and semi-algebraic constraints, these optimization problems can be solved using a matlab based software package named Gloptipoly 3 (see \cite{henrion2003gloptipoly}). Gloptipoly 3 converts the optimization problem into an
equivalent generalized moment problem (GMP) and then computes the global optimal solution(s) by solving a hierarchy of semidefinite program (SDP) relaxations (see \cite{lasserre2001global}, \cite{lasserre2010moments}, \cite{lasserre2006convergent} for more details). It is important to note here, that Gloptipoly solver introduces additional variables as relaxation order increases. If the number of variables in the original optimization problem (e.g., OP1) are $\eta$ and the relaxation order is $\rho$, with $\pi$ number of inequality constraints, the SDP will have $\pi$ semidefinite constraints with moment matrices of size $\omega$ $\times$ $\omega$ where $\omega=$ $\eta+\rho \choose \rho$. Worst case complexity for obtaining an $\epsilon$-optimal solution to a SDP with constraints of size $\omega$ is $O(\sqrt{\omega}\log(\frac{1}{\epsilon}))$ \cite{sturm}. Since the growth of $\omega$ with relaxation order $\rho$ is very fast, using present day desktop computers only small examples can be worked out.  At this point, note that one can also use any standard NLP solvers such as \emph{fmincon}, SNOPT, IPOPT, etc., but at a loss of guarantee of globally optimal solutions.

\subsubsection{Restrictions on number of switchings}
In addition, the computation of time-fuel optimal control with the number of switching or discontinuities restricted as $r\leq n$ can also be handled in the proposed formulation. In such case, we consider sequences for $u^{*}(t)$ from $S_{p,q}$ where $p,q$ are the such that $p,q\leq n$ and $p+q\leq r$. For example, for computing the time-fuel optimal control for a system of order $n=3$ with at most $r=4$ switchings, we use sequences from set $S_{1,3}$, $S_{3,1}$, $S_{2,2}$, $S_{1,2}$ and $S_{2,1}$ and solve ten optimization problem formulated using $u(t)$ as: (1) $\left (+1,0,+1,0,-1 \right )$, (2) $\left (0,+1,0,-1 \right)$ (3) $\left (+1,0,-1,0,+1 \right )$, (4) $\left (+1,0,-1,0,-1 \right)$ and (5) $\left (0,-1,0,+1 \right)$ and their conjugates.

\subsection{Example}
\begin{figure}[h]
\centering
\includegraphics[height=0.3\textwidth, width=0.45\textwidth]{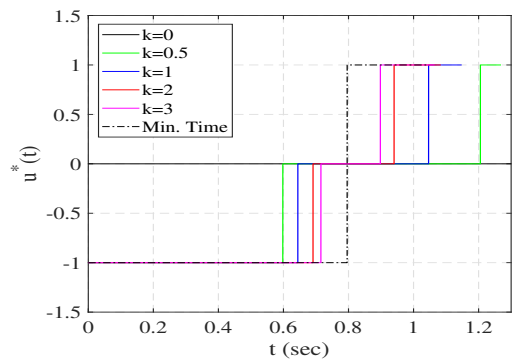} 
\caption{$u^{*}(t)$ Trajectory}
\label{figure3}
\end{figure}

Let us consider a second order system LTI system with $A=$ diag$(-1,-2)$, $B=[1,1]^{T}$ and
We set the initial and final states as: $\mathbf{x}(0)=[0.6,0.4]^{T},~\mathbf{x}(t_f)=\mathbf{0}$. The optimal control $u(t)\in \mathcal{U}^{+}_{c}=\left\{(0,+1,0,-1),(+1,0,+1)\right\}$. Therefore, we achieve two optimizations problems both for problem (OP1) and (OP2), one with $u(t)\in \mathcal{U}^{+}_{c}$ and the other with $u(t)\in \mathcal{U}^{-}_{c}$. 
Similarly, two other optimization problems can be formulated with $u(t)\in \mathcal{U}_c^{-}$. By solving these four optimization problems, the problem (OP1) with $u(t)\in \mathcal{U}_c^{-}$ gives minimum cost and is shown in Figure \ref{figure3}. The performance measures for $k=0,0.5,1,2,3$ and minimum time control are shown in Table \ref{tab:k}. The sparsity in $u^*(t)$ is computed as the ratio of off-duration of $u^*(t)$ to $t_f$. Under the application of $u^*(t)$, the state trajectory steers from $\mathbf{x}_0$ to origin and $u^{*}(t)$ follows one of derived candidate sequence $-1,0,+1$ as shown in Figure \ref{figure3}. 
\begin{table}[h!]
\caption{Performance measures for different values of $k$ in Example}
    \label{tab:k}
    \centering
    \begin{tabular}{c c c c c}
    \hline
       \multirow{2}{*}{$k$}  & \multirow{2}{*}{$J^*$} & \multirow{2}{*}{$t_{f}$} & Time duration & \multirow{2}{*}{Sparsity} \\
       & & & for $u^*(t)\neq 0$ & \\
       \hline 
       0 & 0 & $\infty$ & 0 & 1 \\
      0.5 & 1.2959 & 1.2689 & 0.6615 & 0.4787 \\
        1 & 1.8940 & 1.1480 & 0.746 & 0.3502 \\
        2 & 3.0025 & 1.0839 & 0.8347 & 0.2299   \\
        3 & 4.0752 & 1.0645 & 0.8817 & 0.1717 \\
        Min. Time & -- & 1.0413 & 1.0413 & 0 \\
        \hline
    \end{tabular}
\end{table}
\begin{figure}[h]
\centering
\includegraphics[height=0.3\textwidth, width=0.45\textwidth]{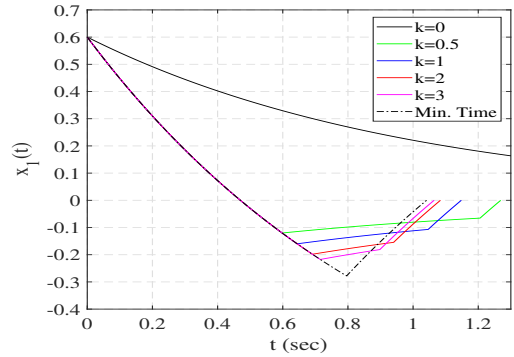} \caption{$x_{1}(t)$ Trajectory}
\label{figure1}
\end{figure}
\begin{figure}[t]
\centering
\includegraphics[height=0.3\textwidth, width=0.45\textwidth]{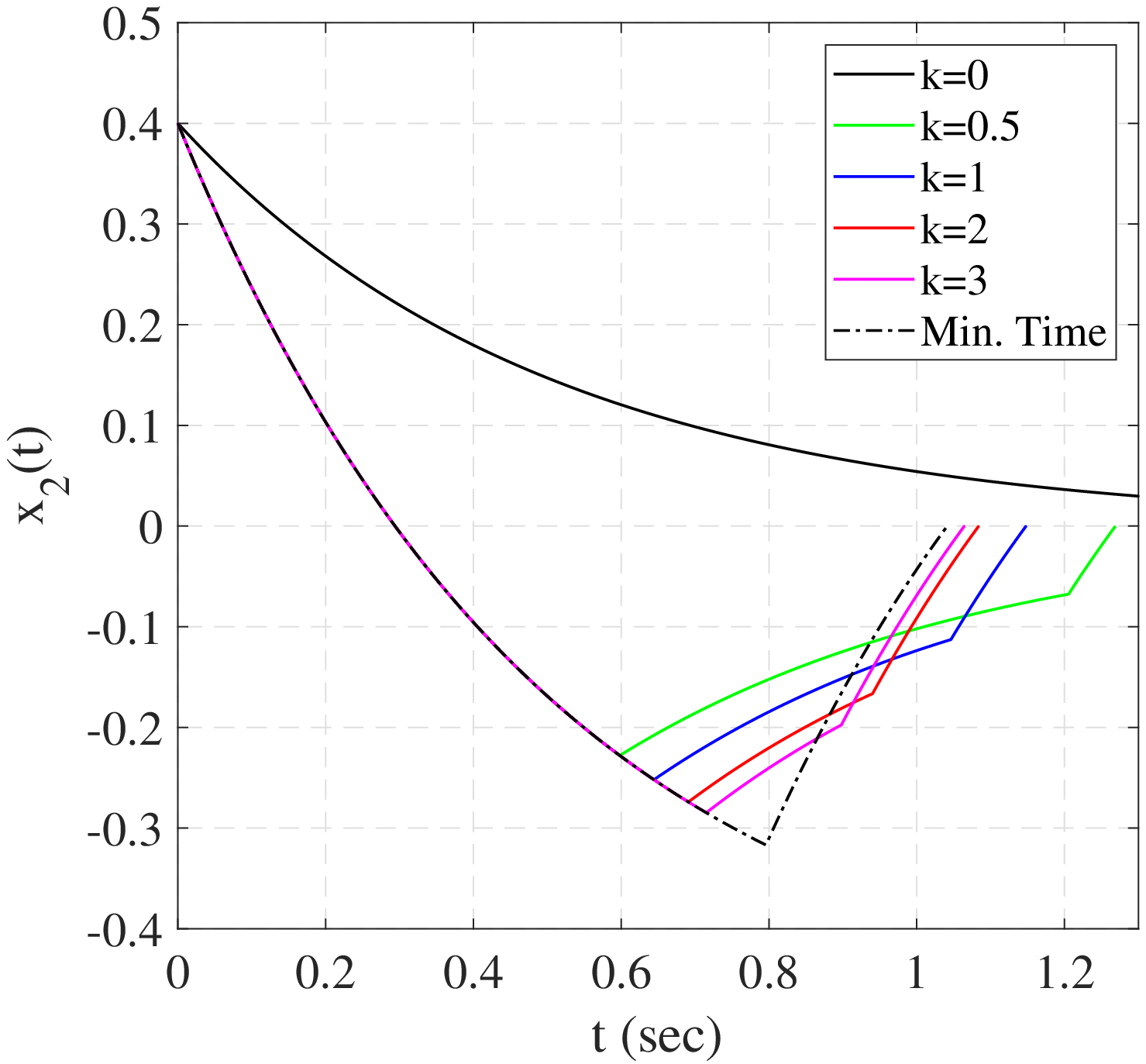} \caption{$x_{2}(t)$ Trajectory}
\label{figure2}
\end{figure}


\section{Conclusion}



 In this article, we computed time-fuel optimal control for LTI systems by characterizing the control in terms of sequences of $+1,0,-1$ and switching time instants.  A method is devised to count and derive all candidate sequences (satisfying PMP necessary conditions). Further, all the candidate sequences are utilized to transform the optimal control problem into multiple static optimization problems which are tractable. Then, the optimal control input is obtained by solving each optimization problem  and selecting the solution with least cost. The computation can be distributed as each optimization problem can be solved separately.  Such characterization of control in terms of time instants can be further exploited in aperiodic feedback control techniques such as self-triggered feedback control \cite{6425820}. 
 
Developing dedicated problem solvers utilizing the structure of cost and constraints is the subject of current and future research. Recently a  way to exploit the sparseness of polynomial constraints  to make this approach scalable for optimal power flow computation appeared in \cite{josz2018lasserre}. Such ideas utilizing any special sparsity structure can be pursued to alleviate the complexity issues the method currently suffers from. Further, a possible classification of initial conditions labelled by the valid candidate sequences is also an interesting direction of research. Such a classification will help in reducing the number of optimization problems that are required to be solved.


\ifCLASSOPTIONcaptionsoff
  \newpage
\fi


%



\addcontentsline{toc}{section}{References}
\bibliographystyle{IEEEtran}
\bibliography{Ref_BIO}
\end{document}